\documentclass[12pt]{amsart}

\newtheorem{theorem}{Theorem}
\newtheorem{theorem*}{Theorem}
\newtheorem{lemma}{Lemma}

\newtheorem{corollary}{Corollary}

\def\mgnb{\overline{\mathcal M}_{g,n}}
\def\mgnpb{\overline{\mathcal M}_{g,n+1}}

\def\mgn{{\mathcal M}_{g,n}}
\def\mgnp{{{\mathcal M}_{g,n+1}}}

\def\mgo{\mathcal{M}_{g,1}}
\def\mgz{\mathcal{M}_{g,0}}
\def\kao{\kappa_1}

\def\mb{\overline{\mathcal M}}

\def\wp{Weil-Pe\-ters\-son}

\begin{document}
\title{Estimates of Weil-Petersson volumes \\ via effective divisors}
\author{Georg Schumacher}
\address{Fachbereich Mathematik und Informatik der Philipps-Universit\"at,
Hans-Meerwein-Strasse, Lahnberge, 35032 Marburg, Germany}
\email{schumac@mathematik.uni-marburg.de}
\author{Stefano Trapani}
\address{Dipartimento di Matematica, Universit\'a di Roma, ``Tor Vergata''
Via della Ricerca Scientifica, 00179 Roma, Italy}
\email{trapani@axp.mat.uniroma2.it}
\maketitle

\section{Introduction and Statement of Results} The total free energy of 
two dimensional gravity, which is a generating function for certain 
intersection numbers on the compactified moduli spaces $\mgnb$ of 
stable $n$-pointed curves, was conjectured by Witten (and proved
by Kontsevich) to satisfy certain KdV equations. This gave new insight 
in the geometry of those moduli spaces \cite{witt, kon, dij}.

The Mumford class $\kao$ on $\mb_{g,0}$ was shown to be proportional to 
the cohomology class $[\omega_{WP}]$ of the Weil-Petersson form by
Wolpert in \cite{wo1}. Furthermore he showed that  the 
restriction of this class to any component of the compactyfying divisor 
coincides with the corresponding Weil-Petersson class. Arbarello and Cornalba 
introduced classes $\kao$ on $\mb_{g,n}$, proved a similar restriction property
for these and concluded proportionality on all  $\mb_{g,n}$ \cite{ac2}. 
So Weil-Petersson volumes ${\rm vol}(\mgn)$ are up to a
normalizing factor the intersection numbers
$$
V_{g, n}=\int_{\mgn}\kao^{3g-3+n}.
$$

Recently, in papers by Kaufmann,
Manin, Zagier, and Zograf \cite{kmz, mz, zo} a generating function was
introduced for intersection numbers of Mumford's tautological classes
\cite{mu2}, and shown to be equal to the above generating function up
to a change of variables.
Previously, Zograf had computed the volumes
for genus 0,1, and 2 explicitly \cite{zo1,zo2}. Manin and Zograf also
gave estimates of the volume for fixed genus and $n\to \infty$.

The aim of this note is to study the \wp\ volume of the moduli spaces
$\mgn$ for fixed $n$ and large $g$. Introducing the decorated
Teichm\"uller space, Penner \cite{pe} gave a technique how to integrate
top degree differential forms on $\mgn$, which led to an estimate of
the volumes of $\mgo$ from below with respect to $ g\to \infty$.
With these methods, Grushevski \cite{gru} recently proved an upper
bound for the volume of $\mgn$ for fixed $n>0$ and large $g$. For $n=1$
his upper estimate has the same order of growth as Penner's lower
estimate.

However, for the classical moduli spaces $\mgz$ the asymptotics of the
volume for $g\to \infty$ have not been treated.

Here, we give a different approach, which does not require the
existence of punctures. On one hand we use the known push-pull type
formulas in the spirit of Arbarello and Cornalba \cite{ac, ac2} to
estimate the volume of $\mgnp$ from below in terms of the volume of
$\mgn$ for any given $g$ and $n \geq 0$. On the other hand, we base
our estimates on the fact that $\kappa_1$ is ample and the above
restriction property.

In this way it is possible to estimate the volume of the moduli space
$\mgz$ from below in terms of the volume of moduli spaces of lower
genus.

We set $V_{(0,3)}=0$. The values $V_{(0,4)}=1$, $V_{(0,5)}=5$, $V_{(1,
1)}=\frac{1}{24}$, and $V_{(1,2)}= \frac{1}{8}$ are known.
We prove the following Theorems.

\begin{theorem}\label{thm1} Let  $2g -2 +n >0$ and $(g,n)\neq
(0,4), (1,1)$. Then
\begin{equation}\label{lower_est}
V_{g,n+1} \geq \frac{1}{2}(3g -2 + n)(7g-7+3n)\cdot V_{g,n} + \frac{1}{24^g g!}
\end{equation}
\end{theorem}

\begin{theorem}\label{thm2} Let $g>1$. Then
\begin{equation}\label{upper_est}
V_{g,0} \geq \frac{1}{28} V_{g-1,2} + \frac{1}{672} V_{g-1,1} +
\frac{1}{14} \sum_{j=2}^{[g/2]} V_{j,1}V_{g-j,1}
-\frac{1}{28} (V_{\frac{g}{2},1})^2 ,
\end{equation}
with $V_{\frac{g}{2},1}=0$, if $g$ is odd.
\end{theorem}

Together with the results of Penner and Grushevsky these imply the
existence of constants $0 < c < C$, independent of $n$ such that
\begin{equation}\label{asymp1}
c^g  (2g)! \leq \frac{V_{g,n}}{(3g-3+n)!} \leq C^g  (2g)!
\end{equation}
for all fixed $n\geq 0$ and large $g$.

In particular for all $n \geq 0$
\begin{equation}\label{asymp2}
\lim_{g \to \infty} \frac{\log  \frac{V_{g,n}}{(3g-3+n)!} }{g \log g} =2
\end{equation}

\section{Proof of the estimates} For any $g$ and $n$ with $2g-2+n>0$ the
map $\mgnp \to \mgn$ forgetting the last puncture is known to extend
holomorphically to a map $$ \pi_{n+1} : \mgnpb \to \mgnb. $$ For $n>0$
it possesses natural sections $\sigma_j$; $j=1,\ldots, n$ (cf.\
\cite[sect.\ 1]{ac2})  with corresponding divisors $D_j$. Denote the
relative dualizing sheaf $\omega_{\mgnpb/\mgnb}$ by $\omega_{n+1}$.
Let
$$
\psi_j:= c_1(\sigma_j^* \omega_{n+1}),
$$
and
$$
K:= c_1(\omega_{n+1}(D)) \in H^2(\mgnpb, {\mathbb R}),
$$
where $D= D_1 + \ldots D_n$.

Finally
$$
\kappa_j:= \pi_{n+1*}(K^{j+1}) \in H^{2j}(\mgnpb, {\mathbb R})
$$
for $j=0,\ldots, 3g-3+n$.

For $n=0$ these are equal to the Mumford classes (denoted also by
$\kappa_j$). Moreover $\kao$ is ample on $\mgnb$.

According to Mumford \cite{mu2} the classes $\kappa_j$ are numerically
effective on $\mb_{g,0}$ for $j=1, \ldots, 3g-3$ in the sense that for
any complete subvariety $W \subset \mb_{g,0}$ of dimension $j$
$$
\int_W \kappa_j \geq 0
$$
holds. Also, for $j_1+\ldots+ j_k = 3g-3$
\begin{equation}\label{harris_int}
\int_{\mb_{g,0}} \kappa_{j_1}\cdot \ldots \cdot \kappa_{j_k} \geq 0
\end{equation}
following \cite{h2}.

We shall need the following formulas for cohomology classes
to be found in
\cite[(1.7), (1.9), and (1.10)]{ac2}.
\begin{gather}
\pi_{n*}(\psi_1^{a_1}\cdot \ldots \cdot \psi_{n-1}^{a_{n-1}}\cdot
\psi_n^{a_n +1})=
\psi_1^{a_1}\cdot \ldots \cdot \psi_{n-1}^{a_{n-1}}\cdot
\kappa_{a_n}  \\
 {\text{\rm for}}\qquad a_j\geq 0\ \nonumber \\
\pi_{n*}(\psi_1^{a_1}\cdot \ldots \cdot \psi_{n-1}^{a_{n-1}})= \hspace{5cm}
\label{2)}
\\
\hspace{2cm}
\sum_{j; a_j>0} \psi^{a_1}_1\cdot \ldots \cdot \psi_{j-1}^{a_{j-1}}\cdot
\psi_j^{a_j -1}\cdot \psi_{j+1}^a\cdot\ldots\cdot \psi_{n-1}^{a_{n-1}}
\nonumber \\
\kappa_a= \pi^*_{n+1}(\kappa_a) + \psi^a_{n+1} \quad \text{\rm on}
\quad \mgnpb \\
\kappa_0 = 2g-2+n. \label{kappa_0}
\end{gather}

\begin{lemma}\label{lem1}
Let $m_j \geq 0$ be integers such that $\sum_{j=1}^k j \cdot m_j = 3g-3
+ n$ with $n\geq 0$. Then
\begin{equation}\label{kappa}
\int_{\mgn} \kappa_1^{m_1}\cdot\ldots\cdot \kappa_k^{m_k} \geq0 .
\end{equation}
\end{lemma}
\begin{proof}The above equation~(\ref{harris_int}) is the statement for
$n=0$. We proceed by induction on $n$. Assume (\ref{kappa}) for some
$n\geq 0$. Then
\begin{gather*}
\int_{\mgnp} \kappa_1^{m_1}\cdot\ldots\cdot \kappa_k^{m_k} =
\hspace{4cm}\\
\int_{\mgn} \pi_{n+1*} \left(\right.
(\pi_{n+1}^* \kappa_1 + \psi_{n+1})^{m_1} \cdot
(\pi_{n+1}^* \kappa_2 + \psi_{n+1}^2)^{m_2} \cdot \\
 \hspace{7cm} \ldots \cdot
(\pi_{n+1}^* \kappa_k + \psi_{n+1}^k)^{m_k}\left.\right)
\end{gather*}
Since
\begin{gather*}
\int_{\mgn} \pi_{n+1*}\left(\pi_{n+1}^*(\kappa_1^{j_1}\cdot
\ldots \cdot \kappa_k^{j_k})\cdot \psi_{n+1}^{\ell + 1}\right) =
\int_{\mgn} \kappa_1^{j_1}\cdot \ldots \cdot \kappa_k^{j_k} \cdot
\kappa_\ell,
\end{gather*}
the above integral can be expressed as a sum of non-negative terms.
\end{proof}
{\it Proof of Theorem~\ref{thm1}.}
We first note that for $g>0$
\begin{gather*}
\int_{\mgnp} \psi_{n+1}^{3g-2+n} = \int_{\mgnp} \psi_{1}^{3g-2+n}
= \int_{\mgn} \psi_{1}^{3g-3+n} = \ldots
= \int_{\mgo} \psi_1^{3g-2}
\end{gather*}
by (\ref{2)}). These integrals are known to be equal to  $1/(24^g
\cdot g!)$ (cf.\ \cite{fp}). For $g=0$
\begin{gather*}
\int_{{\mathcal M}_{0,n+1}} \psi_{n+1}^{n-2} =
\int_{{\mathcal M}_{0,4}}\psi_{1} = \int_{{\mathcal M}_{0,3}}\kappa_0=1 .
\end{gather*}
By Lemma~\ref{lem1}, for $j<3g-2+n$
$$
\int_{\mgnp} \pi_{n+1}^*(\kappa_1)^j\cdot \psi_{n+1}^{3g-2+n-j}=
\int_{\mgn}\kappa_1^j \cdot \kappa_{3g-3+n-j} \geq 0,
$$
hence
\begin{gather*}
\int_{\mgnp}\kappa_1^{3g-2+n} = \int_{\mgnp}\left( \pi_{n+1}^*(\kappa_1)+
\psi_{n+1} \right)^{3g-2+n}= \\
\sum_{j=0}^{3g-2+n} \binom{3g-2+n}{j} \int_{\mgnp} (\pi_{n+1}^*(\kappa_1))^j
\cdot \psi_{n+1}^{3g-2+n-j} \geq \\
(3g-2+n)\cdot \int_{\mgn} \kappa_1^{3g-3+n} \cdot \kappa_0 +
\hspace{4cm}
\\
\frac{1}{2}(3g-2+n)(3g-3+n) \int_{\mgn} \kappa_1^{3g-4+n} \kappa_1 +
\\
\hspace{7cm}
\int_{\mgnp}\psi_{n+1}^{3g-2+n} =\\
((3g-2+n)(2g-2+n) +\frac{1}{2}(3g-2+n)(3g-3+n))\cdot
\int_{\mgn} \kappa_1^{3g-2+n} + \\
\hspace{8cm}
\int_{\mgnp}\psi_{n+1}^{3g-2+n} = \\
\frac{1}{2}(3g-2+n)(7g-7 + 3n) \int_{\mgn} \kappa_1^{3g-3+n} +
\frac{1}{24^g g!}.
\end{gather*}.
\qed

For any family $f:{\mathcal C} \to S$ of stable curves
$\det(f_*\omega_{{\mathcal C}/S})$ is a line bundle over $S$, where
$\omega_{{\mathcal C}/S}$ denotes the relative dualizing sheaf. These
determinant sheaves give rise to a $\mathbb Q$-divisor on $\mb_{g,0}$,
which is usually denoted by $\lambda$. In a similar way the singular
fibers of the above family define devisors, which also give rise to
$\mathbb Q$-divisors on $\mb_{g,0}$, usually denoted by $\delta_i$.
The irreducible components of
the divisor at infinity on $\mb_{g,0}$ are denoted by $\Delta_i$, $i=0,
\ldots, [g/2]$ with classes
$$
[\Delta_i]=\delta_i \ {\text for} \ i\neq1 \ {\text and} \
[\Delta_1]=2 \delta_1.
$$
These are characterized as follows:
\begin{itemize}
\item[(i)]
The generic point of $\Delta_0$ corresponds to an irreducible, stable
curve of genus $g-1$ with one ordinary double point. In fact there is a
generically 2:1 surjective holomorphic map $\mb_{g-1,2} \to \Delta_0$.

\item[(ii)] For $i=1, \ldots, [g/2]$ the generic points of $\Delta_i$
correspond to stable curves with one ordinary double point and two
irreducible components of genus $i$ and $g-i$ resp. There exists a
surjective holomorphic map
$$
\mb_{i,1}\times \mb_{g-i,1} \to \Delta_i,
$$
which is generically 1:1 for $i \neq g/2$ and 2:1 for $i=g/2$
\end{itemize}
(cf.\ \cite{hm}).

\begin{theorem}\label{thm3}
Let
$$
D= p\cdot \lambda - \sum_{j=0}^{[g/2]} q_j\cdot\delta_j \
; \
p, q_j >0
$$
be an effective ${\mathbb Q}$-divisor on $\mb_{g,0}$
such that
$$
\mu_j = \frac{12 q_j -p}{p} >0.
$$
 Then
$$
V_{g,0} > \frac{\mu_0}{2}\cdot V_{g-1,2} + \frac{\mu_1}{48} \cdot V_{g-1,1}
+ \sum_{j=2}^{[g/2]} \mu_j \cdot V_{j,1} \cdot V_{g-j,1}
- \mu_{g/2} \cdot (V_{g/2,1})^2
$$
(where $\mu_{g/2}=V_{g/2,1}=0$, if $g$ is odd).
\end{theorem}

\begin{proof}  According to \cite{mu1}
$$
\kappa_1 = 12 \lambda - \sum_{j=0}^{[g/2]} \delta_j .
$$
We want to write the divisor in the form
$$
D= \alpha \kappa_1 - \sum \beta_j\delta_j
$$
for $\alpha, \beta_j >0$. This gives
$$
\alpha=\frac{p}{12}; \ \beta_j= q_j - \frac{p}{12} >0.
$$
Now
$$
0 < \kappa_1^{3g-4}\cdot D = \alpha \kappa_1^{3g-3} -
\sum \beta_j \kappa_1^{3g-4}\cdot \delta_j.
$$

We use the restriction property of $\kao$ with respect
to $ \Delta_j$.

On $\mb_{g-1,2}$ the class $\kappa_1$ is invariant under the action of
$\mathbb Z_2$ exchanging the punctures. So under the natural map
$\mb_{g-1,2} \to \Delta_0$ it descends to the restriction of the 
class $\kao$ on the ambient space $\mb_{g,0}$. Now
$$
\kappa_1^{3g-4}\cdot \delta_0 = \int_{\Delta_0} \kappa_1^{3g-4}
=\frac{1}{2} V_{g-1,2}.
$$
In a similar way, we get
$$
\kappa_1^{3g-4}\cdot \delta_1 =\frac{1}{2} V_{1,1}\cdot V_{g-1,1},
$$
and for $i>1$
$$
\kappa_1^{3g-4}\cdot \delta_i = V_{i,1}\cdot V_{g-i,1}
$$
with an extra factor $1/2$ for $(V_{\frac{g}{2},1})^2$, if $g$ is
even. Also we use $V_{1,1} = 1/24 $.
\end{proof}

{\it Remark:}\/ Combining push-pull formulas and computations of
 intersections of powers of $\kappa_1$ with various effective
 divisors, one can also estimate the intersection numbers
 $V_{g,n}$ from above in terms of  $V_{g-1,n+2}$, and the numbers
 $V_{j,\ell}$ with $j=g$, $\ell < n$  or  $j<g$, $\ell \leq n+1$. 

The above Theorem~\ref{thm2} follows, since for every rational
$\varepsilon>0$ the $\mathbb Q$-divisor
$D = (11.2 + \varepsilon)\cdot \lambda - \delta$ is ample \cite{mu1},
where $\delta=\sum \delta_j$.

A stronger estimate for $g\geq 23$ follows from the fact that
$\mb_{g,0}$ has positive Kodaira dimension according to the theorems
of Eisenbud, Harris, and  Mumford \cite{eh,h,hm}. The equality
$$
K_{\mb_{g,0}}= 13 \lambda - 2\delta_0 -3\delta_1
-2\sum_{j=2}^{[g/2]} \delta_j
$$
in ${\rm Pic}(\mb_{g,0})\otimes\mathbb Q$ was proved in \cite{hm}.

This implies for $g\geq 23$
$$
V_{g,0} > \frac{11}{26}\cdot V_{g-1,2} + \frac{23}{624} \cdot V_{g-1,1}
+ \frac{11}{13}
\sum_{j=2}^{[g/2]}  \cdot V_{j,1} \cdot V_{g-j,1}
- \frac{11}{26} \cdot (V_{g/2,1})^2.
$$

Further computations of effective divisors in terms of $\lambda$ and $\delta_j$ 
were provided in \cite{hm,eh}. All of these divisors
satisfy the hypothesis of Theorem~\ref{thm3} and can be
used for new estimates of the Weil-Petersson volumes.

We finally discuss, how to arrive at the asymptotic estimates
(\ref{asymp1}).

A rough estimate following from Theorem~\ref{thm1} is
$$
V_{g,n} \geq V_{g,1},
$$
which, together with Penner's lower estimate for $V_{g,1}$,
already implies the existence of a constant $c>0$, independent of
$n\geq 1$, such that for large $g$
$$
\frac{V_{g,n}}{(3g-3+n)!} \geq c^g (2g)!.
$$
The corresponding upper estimate is due to Grushevski,
so that (\ref{asymp1}) follows for $n\geq 1$.

For $n=0$ and $g >1$ the above
Theorem~\ref{thm1} and Theorem~\ref{thm2} give
$$
\frac{1}{672} V_{g-1,1} \leq V_{g,0} < \frac{2}{(3g-2)(7g-7)} V_{g,1}.
$$
Again, with \cite{gru,pe} these inequalities yield
the following
corollary.
\begin{corollary}
There exist constants $0< \tilde c <\tilde C$ such that for $g \gg 0$
$$
\tilde c^{g} (2g)! \leq \frac{V_{g,0}}{(3g-3)!} \leq \tilde C^{g} (2g)!.
$$
\end{corollary}

\begin{small}
This paper was written, while the second named author was visiting the 
University of Marburg. He would like to thank the Department of Mathematics 
for its hospitality. 
\end{small}


\begin{thebibliography}{MMMM}

\bibitem[A-C 1]{ac} Arbarello, E., Cornalba, M.:
The Picard groups of the moduli spaces of curves.
Topology {\bf 26} (1987) 153--171.

\bibitem[A-C 2]{ac2} Arbarello, E., Cornalba, M.: Combinatorial and
algebro geometric cohomology classes on the moduli space of curves. J.\
Alg.\ Geom.\ {\bf 5} (1996) 705--749.

\bibitem[DIJ]{dij} Dijkgraaf, R.: Intersection theory, integrable
hierarchies and topological field theory. In; new symmetric principles
in quantum field theory, Carg\'ese, 1991, Adv.\ Sci.\ Int.\ Ser.\ B
Phys.\ {\bf 295}, Plenum, New York, 1992.

\bibitem[E-H]{eh} Eisenbud, D., Harris, J.:
The Kodaira dimension of the moduli space of curves of genus $\ge 23$.
Invent.\ Math.\ {\bf 90} (1987) 359--387.

\bibitem[F-P]{fp} Faber, C., Pandharipande, R.:
Hodge integrals and Gromov-Witten theory.
Invent.\ Math.\ {\bf 139} (2000) 173-199.

\bibitem[GE]{getz} Getzler, E.: Intersection theory on
$\overline{\mathcal{M}}_{1,4}$ and elliptic Gromov-Witten invariants.
J.\ Am.\ Math.\ Soc.\ {\bf 10} (1997) 973-998.

\bibitem[GR]{gru} Grushevsky, S.:
Explicit upper bound for the Weil-Petersson volumes. Preprint
math.AG/0003217

\bibitem[HA1]{h} Harris, J.:
On the Kodaira dimension of the moduli space of curves. II:
The even-genus case.
Invent.\ Math.\ {\bf 75} (1984) 437-466.

\bibitem[HA2]{h2} Harris, J.: Families of smooth curves.
Duke Math.\ J.\ {\bf 51} (1984) 409-419.


\bibitem[H-M]{hm} Harris, J., Mumford, D.:
On the Kodaira dimension of the moduli space of curves.
Invent.\ Math.\ {\bf 67} (1982) 23--86.

\bibitem[K-M-Z]{kmz} Kaufmann, R.; Manin, Yu.; Zagier, D.:
Higher Weil-Petersson volumes of moduli spaces of stable $n$-pointed curves.
Commun.\ Math.\ Phys.\ {\bf 181} (1996) 763--787.

\bibitem[KO]{kon} Kontsevich , M.: Intersection theory on the moduli
space of curves and the matrix Airy function. Comm.\ Math. Phys. {\bf
147} (1992) 1--23.

\bibitem[M-Z]{mz} Manin, Y., Zograf,P.:
Invertible Cohomological Field Theories and Weil-Petersson volumes.
Preprint math.AG/9902051.

\bibitem[MU1]{mu1} Mumford, D.:
Stability of projective varieties.
Enseign. Math., II.\ Ser.\ {\bf 23} (1977) 39--110.


\bibitem[MU2]{mu2} Mumford, D.:
Towards an enumerative geometry of the moduli space of curves.
Arithmetic and geometry, Pap.\ dedic.\ I.R.\ Shafarevich, Vol. II:
Geometry, Prog.\ Math.\ {\bf 36} (1983) 271--328.

\bibitem[PE]{pe} Penner, R.C.: \wp\ volumes. J.\ Diff.\ Geom.\ {\bf
35} (1992) 559--608.


\bibitem[WI]{witt} Witten, E.: Two dimensional gravity and intersection
theory on moduli spaces. Surveys in Diff.\ Geom.\ {\bf 1}  (1991)
243--310.

\bibitem[WO]{wo1} Wolpert, S.: On the homology of the moduli spaces of
stable curves. Ann.\ Math.\ {\bf 118} (1983) 491--523.


\bibitem[ZO1]{zo1} Zograf, P.: The Weil-Petersson volume of the moduli
space of punctured spheres. B\"odigheimer, et al.,
Mapping class groups and moduli spaces of Riemann surfaces.
Proceedings of workshops G\"ottingen 1991, Seattle, 1991, American
Math.\ Soc., Contemp.\ Math.\ {\bf 150} (1993) 367--372.


\bibitem[ZO2]{zo2} Zograf, P.: Weil-Petersson volumes of low genus
moduli spaces, Funct.\ Anal.\ Appl.\ {\bf 32} (1998) 78--81.

\bibitem[ZO3]{zo} Zograf, P.: Weil-Petersson volumes of moduli spaces of
curves and the genus expansion in two dimensional gravity. Preprint
math.AG/9811026.

\end{thebibliography}
\end{document}